\documentclass{article}
\usepackage{graphicx, subfigure, amsfonts, amsmath, amsthm, tikz, url}
\usepackage{a4, color}
\usepackage[top=1in, bottom=1.25in, left=1.25in, right=1.25in]{geometry}
\newcommand{\ignore}[1]{}
\newcommand{\paf}[2]{\frac{\partial #1}{\partial #2}}
\newcommand{\naf}[2]{\frac{\text d #1}{\text d #2}}

\newcommand{\diag}{\mbox{diag}}

\newcommand{\ContinuousDiscrete}[2]{
    ~\\[1ex]
    \begin{tabular}{c|c}
    \multicolumn{1}{c}{\bf Continuous} & \multicolumn{1}{c}{\bf Discrete} \\[-2ex]
    \begin{minipage}{7.2cm} 
    {\small #1}\\[-4ex]
    \end{minipage}
    &
    \begin{minipage}{7.2cm} 
    {\small #2}\\[-4ex]
    \end{minipage}
    \end{tabular}\\[1ex]
}

\begin{document}
\title{Derivations of continuous and discrete energy equations in wave and shallow-water equations}
\author{Bas van 't Hof\footnote{Corresponding author. bas.vanthof@vortech.nl} \qquad Mathea~J. Vuik\footnote{thea.vuik@vortech.nl. VORtech, Westlandseweg 40d, 2624AD Delft, The Netherlands.}}

\maketitle

\begin{abstract}
  Symmetry-preserving (mimetic) discretization aims to preserve
  certain properties of a continuous differential operator in its
  discrete counterpart.  For these discretizations, stability and
  (discrete) conservation of mass, momentum and energy are proven
  in the same way as for the original continuous model.  
  
  In our papers \cite{Hof17V} and \cite{Hof19V}, we presented space discretization schemes for 
  various models, which had exact conservation of mass, momentum and energy.
  Mass and momentum conservation followed from the left null spaces of the discrete operators used.
  The conservation of energy in the continuous and discrete models is more complicated, 
  and the papers had little space for their complete derivation. This paper contains the 
  derivation of the energy equations in more detail than was given in
  the papers \cite{Hof17V} and \cite{Hof19V}.
  \end{abstract}

\textbf{Symmetry-preserving discretizations, Mimetic methods, Finite-difference methods, Mass, momentum and energy conservation, Curvilinear staggered grid}

\maketitle

\section{Introduction and motivation}
All of the models presented in \cite{Hof17V} and \cite{Hof19V}, except the scalar wave equation, consist of continuity, momentum and state equations. 
The energy equation is derived by combining these three equations. All the continuous energy equations express the change in the energy density $e$ 
in terms of energy fluxes $\vec f_e$, and therefore have the general form
\begin{equation}
\paf et + \nabla \cdot \vec f_e = 0. 
\label{eq: local energy balance}
\end{equation}
The time derivative of the total energy $E$, which is the integral of the energy density over a domain $V$, has only boundary terms:
\[ \paf Et + \oint_{\delta V}~\vec f_e \mbox dS = 0, \]
where $\delta V$ is the boundary of the domain $V$. In cases without boundary effects, such as periodic domains and
domains with energy-conserving boundary conditions, the total energy remains constant.

In the discrete models, the energy is located on all the points in the staggered grid, and a local energy balance like (\ref{eq: local energy balance})
cannot be given. Instead, the conservation of total discrete energy ${\tt E}$ is shown by deriving 
\begin{equation} 
\paf {\tt E}t = 0,
\label{eq: discrete energy equation}
\end{equation}
using certain properties of the discrete operators used in the discretizations, like the discrete Laplacian {\sf LAPL}, the discrete divergence {\sf DIV}, 
the discrete gradient {\sf GRAD} and others.

The subsequent four sections each present one of the models. In each case a continuous and a discrete model is presented, and the energy equations are derived.

\section{Energy equation in scalar wave equations}
The scalar wave equation describes the change in the pressure $p$ or its discrete approximation {\tt p}, and is given by
\ContinuousDiscrete{
\begin{eqnarray}
 \paf{{}^2p}{t^2} = \nabla^2 p.
\end{eqnarray}
}{
\begin{eqnarray}
 \naf{{}^2{\tt p}}{t^2} = {\sf LAPL}~{\tt p},
\end{eqnarray}
}
where ${\sf LAPL}$ is the discrete approximation of the Laplacian operator $\nabla^2$, which is symmetrical:
\[
{\sf LAPL}^* = {\sf LAPL}.
\]

In the continuous model, the total energy $E$ is the integral of the energy density $e$. 
In the discrete model, the total energy $E$ is presented using scalar products:
\ContinuousDiscrete{
\begin{eqnarray}
 e := \frac 12 \left( \paf pt \right)^2 + \frac 12 |\nabla p|^2
\end{eqnarray}
}{
\begin{eqnarray}
 E := \frac 12 \left< \naf {\tt p}t, \naf {\tt p}t \right>_c - \frac 12 \left<{\tt p}, {\tt LAPL}~{\tt p}\right>_c.
\end{eqnarray}
}
The time derivative of the energy is
\ContinuousDiscrete{
\begin{eqnarray}
\paf et &=& 
\paf pt \paf{{}^2p}{t^2} + \nabla p \cdot \paf{}t \nabla p  
\\ &=&
\nabla \cdot \left( \paf pt \nabla p\right). \nonumber 
\end{eqnarray}~\\[3ex]
}{
\begin{eqnarray}
\naf Et &=& 
  \left< \naf {\tt p}t, \naf {{}^2 \tt p}{t^2} \right>_c \\ &&
  - \frac 12 \left<\naf {\tt p}t, {\sf LAPL}~{\tt p}\right>_c
  - \frac 12 \left<{\tt p}, {\sf LAPL}~\naf {\tt p}t \right>_c
\nonumber \\ &=&
  \left< \naf {\tt p}t, {\sf LAPL}~{\tt p} \right>_c 
  - \frac 12 \left<\naf {\tt p}t, ({\sf LAPL}+{\sf LAPL}^*) {\tt p}\right>_c.\nonumber
\end{eqnarray}
}
Using the symmetry property that ${\sf LAPL}^*={\sf LAPL}$, the following energy equation is found:
\ContinuousDiscrete{
\begin{eqnarray}
\paf et &+& 
\nabla \cdot \left(-\paf pt \nabla p\right) = 0.
\end{eqnarray}
}{
\begin{eqnarray}
\paf Et=0.
\end{eqnarray}
}

\section{Energy equation in linear-wave equations}
The linear-wave equations describe the change in the flow velocity $\vec v$, the density $\rho$ and the pressure $p$, 
and their discrete approximations {\tt v}, {\tt rho} and {\tt p}. The equations are given in the form of the continuity, momentum and state equations
\ContinuousDiscrete{
\begin{eqnarray}
\paf \rho t &+& \nabla \cdot \rho_0 \vec v = 0 \nonumber \\ 
\paf {\rho_0 \vec v}t &+& \nabla p= 0,  \nonumber \\
\rho &=& c^2 p.
\end{eqnarray}
}{
\begin{eqnarray}
\paf{}t{\tt rho}  &+& \rho_0 {\sf DIV}~{\tt v} = 0 \nonumber \\ 
\rho_0 \paf{\tt v}t &+& {\sf GRAD}~{\tt p}= 0,  \nonumber \\
{\tt rho} &=& c^2 {\tt p}.
\end{eqnarray}
}
where $c$ is the wave propagation speed, $\rho_0$ is a constant reference density, ${\sf DIV}$ is the discrete approximation of the divergence $\nabla \cdot$, and {\sf GRAD} of the gradient $\nabla$.
The discrete divergence and gradient are each other's negative adjoint:
\[ {\sf GRAD}^* = -{\sf DIV}.\]

\clearpage
\noindent{\bf Time derivative of kinetic energy}\\[1ex]
The local kinetic energy $e_{kin}$ and the total kinetic energy $E_{kin}$ are given by
\ContinuousDiscrete{
\begin{eqnarray}
e_{kin} := \frac{\rho_0} 2 |\vec v|^2.
\end{eqnarray}
}{
\begin{eqnarray}
E_{kin} := \frac{\rho_0} 2 \left<{\tt v},{\tt v}\right>_v.
\end{eqnarray}
}
The time derivative of the kinetic energy is given by
\ContinuousDiscrete{
\begin{eqnarray}
\paf{e_{kin}}t &=& \rho_0 \vec v \cdot \paf{\vec v}t .
\label{d ekin dt lw}
\end{eqnarray}
}{
\begin{eqnarray}
\paf{E_{kin}}t &=& \rho_0 \left<{\tt v},\naf {\tt v}t \right>_v.
\label{d ekin dt lw discrete}
\end{eqnarray}
}
The time derivatives in the right-hand sides of (\ref{d ekin dt lw}-\ref{d ekin dt lw discrete}) are eliminated using the momentum equation
and the following expression is found
\ContinuousDiscrete{
\begin{eqnarray}
\paf{e_{kin}}t &=& -\vec v \cdot \nabla p. 
\end{eqnarray}
}{
\begin{eqnarray}
\paf{E_{kin}}t &=&  -\left< {\tt v},{\sf GRAD}~{\tt p}\right>_v
\end{eqnarray}
}
{\bf Internal energy}\\[1ex]
The local internal energy $e_{int}$ and total internal energy $E_{int}$ are defined by
\ContinuousDiscrete{
\begin{eqnarray}
 e_{int} := \frac{c^2}{2\rho_0}\rho^2,
\end{eqnarray}
}{
\begin{eqnarray}
 E_{int} := \frac{c^2}{2\rho_0}\left< {\tt rho} ,{\tt rho} \right>_c.
\end{eqnarray}
}
Their time derivatives are given by
\ContinuousDiscrete{
\begin{eqnarray}
 \paf{}t e_{int} = \frac{c^2}{\rho_0}\rho \paf \rho t,
\end{eqnarray}
}{
\begin{eqnarray}
 \naf{}t E_{int} := \frac{c^2}{\rho_0}\left< {\tt rho} ,\naf{}t {\tt rho} \right>_c.
\end{eqnarray}
}
The time derivatives in the right-hand sides are eliminated using the continuity equation:
\ContinuousDiscrete{
\begin{eqnarray}
 \paf{}t e_{int} = -c^2\rho \nabla \cdot \vec v = 
                   -p \nabla \cdot \vec v,
\end{eqnarray}
}{
\begin{eqnarray}
 \naf{}t E_{int} := -\left< {\tt p} ,{\sf DIV}~{\tt v} \right>_c.
\end{eqnarray}
}
{\bf Energy equation}\\[1ex]
The local energy $e=e_{kin}+e_{int}$ is the sum of local kinetic and internal energies, and the total energy $E=E_{kin}+E_{int}$ is the sum of the 
total kinetic and internal energies, so their time derivatives are
\ContinuousDiscrete{
\begin{eqnarray}
 \paf et  = -\vec v \cdot \nabla p -p \nabla \cdot \vec v = -\nabla (p \vec v).
\end{eqnarray}
}{
\begin{eqnarray}
 \naf Et &=&
-\left< {\tt v},{\sf GRAD}~{\tt p}\right>_v
-\left< {\tt p} ,{\sf DIV}~{\tt v} \right>_c \nonumber \\ &=&
- \left< {\tt p} ,({\sf GRAD}^*+{\sf DIV})~{\tt v} \right>_c.
\end{eqnarray}
}
Using the symmetry property ${\sf GRAD}^*=-{\sf DIV}$, the following energy equation is found:
\ContinuousDiscrete{
\begin{eqnarray}
    \paf et &+& \nabla \cdot p \vec v  =0.
\end{eqnarray}
}{
\begin{eqnarray}
 \naf Et = 0.
\end{eqnarray}
}

\clearpage
\section{Energy equation in compressible-wave equations}
\noindent
{\bf Compressible-wave equations}\\[1ex]
The compressible-wave equations are given by the continuity, momentum and state equations
\ContinuousDiscrete{
\begin{eqnarray}
\paf \rho t &+& \nabla \cdot \rho \vec v = 0, \nonumber \\
\paf {\vec v}t &+& \nabla Q(p) = 0, \nonumber \\
\rho &=& R(p).
\end{eqnarray}
}{
\begin{eqnarray}
\naf{}t {\tt rho} &+& {\sf DIV}\tilde{\sf r} ~{\tt v} = 0, \nonumber \\
\naf{\tt v}t &+& {\sf GRAD}~Q({\tt p}) = 0, \nonumber \\
{\tt rho} &=& R({\tt p}).
\end{eqnarray}
}
where the function $Q$ is given in terms of the density function $R$ as $Q(p) := \int^p \frac 1{R(p)}~\mbox dq$, so the momentum equation may also be written as
\begin{eqnarray}
\paf {\vec v}t + \frac 1\rho \nabla p= 0.
\end{eqnarray}
Another form of the momentum equation uses the function $S(p) := \int^p \frac 1{R^2(p)}~\mbox dq$, and reads
\ContinuousDiscrete{
\begin{eqnarray}
\paf {\vec v}t + \rho \nabla S(p)= 0, & 
\end{eqnarray}
}{
\begin{eqnarray}
\naf{\tt v}t &+& \tilde{\sf r}{\sf GRAD}~S({\tt p}) = 0.
\end{eqnarray}
}
The operator $\tilde{\sf r}{\sf GRAD}$, which is the discrete approximation of the operator $\rho \nabla$, is related to the discrete gradient ${\sf GRAD}$ 
in the discrete chain rule
\[
\tilde{\sf r}{\sf GRAD}~S({\tt p}) = {\sf GRAD}~Q({\tt p}),
\]
and the operator ${\sf DIV}\tilde{\sf r}$, the discrete approximation of the operator $\nabla \cdot \rho$, is given by
\[
{\sf DIV}\tilde{\sf r} = -\tilde{\sf r}{\sf GRAD}^*.
\]

\noindent
{\bf Time derivative of kinetic energy}\\[1ex]
The local kinetic energy $e_{kin}$ and the total kinetic energy $E_{kin}$ are given by
\ContinuousDiscrete{
\begin{eqnarray}
e_{kin} := \frac {\rho_0} 2 |\vec v|^2,
\end{eqnarray}
}{
\begin{eqnarray}
E_{kin} := \frac {\rho_0} 2 \langle {\tt v}, {\tt v} \rangle_v,
\end{eqnarray}
}
The time derivative of the kinetic energy is given by
\ContinuousDiscrete{
\begin{eqnarray}
\paf{e_{kin}}t &=& \rho_0 \vec v \cdot \paf{\vec v}t.
\label{d ekin dt cw}
\end{eqnarray}
}{
\begin{eqnarray}
\paf{E_{kin}}t &=& \rho_0 \left< {\tt v}, \naf {\tt v}t \right>_v.
\label{d ekin dt cw discrete}
\end{eqnarray}
}
{\bf Time derivative of kinetic energy converted to spatial derivatives}\\[1ex]
The time derivatives in the right-hand sides of (\ref{d ekin dt cw}-\ref{d ekin dt cw discrete}) are replaced by the expression given in the momentum equation
and the following expression is found
\ContinuousDiscrete{
\begin{eqnarray}
\paf{e_{kin}}t &=& 
- \rho_0 \vec v \cdot \nabla Q(p).
\end{eqnarray}
}{
\begin{eqnarray}
\paf{E_{kin}}t &=& - \rho_0 \left< {\tt v}, {\sf GRAD}~Q({\tt p}) \right>_v.
\end{eqnarray}
}
{\bf Internal energy}\\[1ex]
The local internal energy $e_{int}$ is given by
\begin{eqnarray}
   e_{int} &:=& \rho_0 \int^p \frac{R(p)-R(q)}{R^2(q)}\mbox dq =
\rho_0 R(p) \int^p \frac 1{R^2(q)} \mbox dq -
\rho_0 \int^p \frac 1{R(q)} \mbox dq. 
\end{eqnarray}
Its derivative with respect to the pressure is given by
\begin{eqnarray}
e'_{int}(p) &=&  
\rho_0 R'(p) \int^p \frac 1{R^2(q)} \mbox dq + \rho_0 R(p) \frac 1{R^2(p)} - \rho_0 \frac 1{R(p)} = 
\rho_0 R'(p) \int^p \frac 1{R^2(q)} \mbox dq
\nonumber \\  &=& \rho_0 R'(p) S(p).
\end{eqnarray}
{\bf Time derivative of the internal energy}\\[1ex]
The chain rule is applied to find the following expression for the time derivative of the internal energy:
\begin{eqnarray}
\paf{e_{int}}t &=&  
e'_{int}(p) \paf pt = 
\rho_0 R'(p) S(p) \paf pt = 
\rho_0 S(p) \paf \rho t.
\end{eqnarray}
Using the continuity equation, the time derivative is eliminated
\ContinuousDiscrete{
\begin{eqnarray}
\paf{e_{int}}t &=&  
- \rho_0 S(p) \nabla \cdot \rho \vec v 
\end{eqnarray}~\\[5ex]
}{
\begin{eqnarray}
\paf{{\tt e}_{int}}t &=&  
-\rho_0~\diag\left(S({\tt p}) \right) ~{\sf DIV}\tilde{\sf r}~{\tt v}.
\nonumber \\
\paf{E_{int}}t &=&  
-\rho_0 \langle {\tt c1}, 
\diag\left(S({\tt p}) \right) ~{\sf DIV}\tilde{\sf r}~{\tt v}\rangle_c
\nonumber \\ &=&
-\rho_0 \langle S({\tt p}), {\sf DIV}\tilde{\sf r}~{\tt v}\rangle_c.
\end{eqnarray}
}
{\bf Energy equation}\\[1ex]
The time derivatives of local and total energies $e$ and $E$ are 
\ContinuousDiscrete{
\begin{eqnarray}
 \paf et &=& 
- \rho_0 \vec v \cdot \nabla Q(p)
- \rho_0 S(p) \nabla \cdot \rho \vec v.
\end{eqnarray}\\[2ex]
}{
\begin{eqnarray}
 \paf Et &=&  
- \rho_0 \left< {\tt v}, {\sf GRAD}~Q({\tt p}) \right>_v
-\rho_0 \langle S({\tt p}), {\sf DIV}\tilde{\sf r}~{\tt v}\rangle_c
\nonumber \\ 
\end{eqnarray}
}
Now we use 
the symmetry property that ${\sf DIV\tilde r}^*=-{\sf \tilde rGRAD}$ 
and 
the chain rules $\nabla Q = \rho \nabla S$, ${\sf GRAD} Q({\tt p}) = 
\tilde{\sf r}{\sf GRAD}~S({\tt p})$,
to find
\ContinuousDiscrete{
\begin{eqnarray}
 \paf et &=& 
- \rho_0 \rho \vec v \cdot \nabla S(p)
- \rho_0 S(p) \nabla \cdot \rho \vec v   
\nonumber \\ &= &
- \rho_0 \nabla \cdot (\rho \vec v~S(p) )
\end{eqnarray}~\\[2ex]
}{
\begin{eqnarray}
 \paf Et &=&  
- \rho_0 \left< {\tt v}, {\sf GRAD}~Q({\tt p}) \right>_v
+\rho_0 \langle \tilde{\sf r}{\sf GRAD}~S({\tt p}), {\tt v}\rangle_v
\nonumber \\ &=&
- \rho_0 \left< {\tt v}, {\sf \tilde rGRAD}~S({\tt p}) \right>_v
+\rho_0 \langle \tilde{\sf r}{\sf GRAD}~S({\tt p}), {\tt v}\rangle_v
\nonumber \\
\end{eqnarray}
}
The energy equation is therefore
\ContinuousDiscrete{
\begin{eqnarray}
 \paf et &+& 
\rho_0 \nabla \cdot \left( \rho \vec v S(p) \right) = 0.
\end{eqnarray}
}{
\begin{eqnarray}
 \paf Et =0. 
\end{eqnarray}
}

\clearpage
\section{Energy equation in isentropic compressible Euler equations}

\noindent
{\bf Isentropic compressible Euler equations}\\[1ex]
The isentropic compressible Euler equations are given by the continuity, momentum and state equations
\ContinuousDiscrete{
\begin{eqnarray}
\paf \rho t &+& \nabla \cdot \rho \vec v = 0, \nonumber \\
\paf {\rho \vec v}t &+& \nabla \cdot \rho \vec v \otimes \vec v + \nabla p= 0,\nonumber \\
\rho &=& R(p).
\end{eqnarray}
}{
\begin{eqnarray}
\naf{}t {\tt rho} &+& {\sf DIVr}~{\tt v} = 0, \nonumber \\ 
\naf{\tt rv}t &+& {\sf ADVEC}~{\tt v} + {\sf GRAD}~{\tt p} = 0, \nonumber \\
{\tt rho} &=& R({\tt p}),
\end{eqnarray}
}
where ${\sf ADVEC}$ is the discrete approximation of the advection operator $\nabla \cdot \rho v \otimes$, 
{\sf DIVr} of the operator $\nabla \cdot \rho$, ${\sf rGRAD}$ of $\rho \nabla$, 
and where the discrete local momentum ${\tt rv}$ is given by
\begin{eqnarray}
 {\tt rv} :=  \diag\left({\tt v}\right) ~{\sf Interp}_{v\leftarrow c} {\tt rho},
\end{eqnarray}
where ${\sf Interp}_{v\leftarrow c}$ is an interpolation which uses the densities at the cell-centers of the staggered grid to calculate densities at the cell-faces.

The operators ${\sf DIVr}$ and ${\sf rGRAD}$ are each other's negative adjoints:
\[ {\sf rGRAD}^* = -{\sf DIVr},\]
and the operator ${\sf rGRAD}$ is related to the discrete gradient ${\sf GRAD}$ in the discrete chain rule
\[ 
{\sf GRAD}~{\tt p}={\sf rGRAD}~Q({\tt p}).\]
The advection operator has the following symmetry property:
$${\sf ADVEC}+{\sf ADVEC}^*=\diag({\sf Interp}_{v\leftarrow c} {\sf DIVr}~{\tt v}).$$

{\bf Time derivative of kinetic energy}\\[1ex]
The local continuous kinetic energy $e_{kin}$ and the total discrete kinetic energy $E_{kin}$ are given by
\ContinuousDiscrete{
\begin{eqnarray}
e_{kin} := \frac \rho 2 |\vec v|^2,
\end{eqnarray}
}{
\begin{eqnarray}
E_{kin} := \frac 12 \left< {\tt v}, {\tt rv} \right>_v.
\end{eqnarray}
}
Using the product rule for differentiation, the time derivative of the 
kinetic energy is given by
\ContinuousDiscrete{
\begin{eqnarray}
\paf{e_{kin}}t &=& \vec v \cdot \paf{\rho \vec v}t - \frac {|\vec v|^2}2 \paf \rho t,
\label{d ekin dt}
\end{eqnarray}\\[2ex]
}{
\begin{eqnarray}
\paf{E_{kin}}t &=& 
\left< {\tt v}, \naf{\tt rv}t \right>_v -  \label{d ekin dt discrete}
\\ && \frac 12
\left< \diag({\tt v}){\tt v}, {\sf Interp}_{v\leftarrow c}\naf{}t  {\tt rho} \right>_v. 
 \nonumber
\end{eqnarray}
}
{\bf Time derivative of kinetic energy converted to spatial derivatives}\\[1ex]
The time derivatives in the right-hand sides of (\ref{d ekin dt}-\ref{d ekin dt discrete}) are eliminated using the continuity and momentum equations
and the following expression is found
\ContinuousDiscrete{
\begin{eqnarray}
\paf{e_{kin}}t &=& 
-\vec v \cdot \nabla \cdot \rho \vec v \otimes \vec v
-\vec v \cdot \nabla p
\nonumber \\ &&
+ \frac {|\vec v|^2}2 \nabla \cdot \rho \vec v,
\end{eqnarray}
}{
\begin{align}
\paf{E_{kin}}t &= 
-\left< {\tt v}, {\sf ADVEC}~{\tt v} \right>_v 
-\left< {\tt v}, {\sf GRAD} ~{\tt p} \right>_v 
\nonumber \\ &
+\frac 12
\left<  \diag({\tt v}){\tt v}, {\sf Interp}_{v\leftarrow c}{\sf DIVr}~{\tt v} \right>_v.
\end{align}
}

To derive the local energy balance, we need the product rule for advection, given by
\begin{eqnarray} 
   \nabla \cdot \rho |\vec v|^2/2 \vec v &=& 
   \vec v \cdot \nabla \cdot \rho \vec v \otimes \vec v - 
   |\vec v|^2 \nabla \cdot \rho \vec v/2.
\end{eqnarray}
Using this rule and the chain rules that $\nabla p = \rho \nabla Q$ and ${\sf GRAD}~{\tt p}={\sf rGRAD}~Q({\tt p})$, it is found that
\ContinuousDiscrete{
\begin{align}
\paf{e_{kin}}t &=
 -\nabla \cdot \frac \rho 2 |\vec v|^2 \vec v 
 -|\vec v|^2 \nabla \cdot \frac \rho 2  \vec v
\nonumber \\ & -
\rho \vec v \cdot \nabla Q(p)
\nonumber \\ & +
\frac {|\vec v|^2}2 \nabla \cdot \rho \vec v,
\end{align}
}{
\begin{align}
\paf{E_{kin}}t &=
-\frac 12 \left< {\tt v}, ({\sf ADVEC}+{\sf ADVEC}^*)~{\tt v} \right>_v 
\nonumber \\ &-
\left< {\tt v}, {\sf rGRAD} ~Q({\tt p}) \right>_v 
\nonumber \\ &+
\frac 12
\left<  \diag({\tt v}){\tt v}, {\sf Interp}_{v\leftarrow c}{\sf DIVr}~{\tt v} \right>_v.
\end{align}
}
The second and last terms in the continuous equation cancel each other. Also, the first and last 
terms in the discrete equation cancel, because of the symmetry property 
$${\sf ADVEC}+{\sf ADVEC}^*=\diag({\sf Interp}_{v\leftarrow c} {\sf DIVr}~{\tt v}).$$
This leads to the shorter equations
\ContinuousDiscrete{
\begin{eqnarray}
\paf{e_{kin}}t &=& 
 -\nabla \cdot \frac \rho 2 |\vec v|^2 \vec v -
\rho \vec v \cdot \nabla Q(p),
\end{eqnarray}
}{
\begin{eqnarray}
\paf{E_{kin}}t &=& -
\left< {\tt v}, {\sf rGRAD} ~Q({\tt p}) \right>_v.
\end{eqnarray}
}
{\bf Time derivative of internal energy}\\[1ex]
The local internal energy $e_{int}$ is given by
\begin{eqnarray}
   e_{int} &=& \int^p \frac{R(p)-R(q)}{R(q)}\mbox dq = 
                R(p) \int^p \frac 1{R(q)}\mbox dq  - p.
\end{eqnarray}
Its derivative with respect to the pressure is given by
\begin{eqnarray}
e'_{int}(p) &=& 
               R'(p) \int^p \frac 1{R(q)}\mbox dq + R(p) \frac 1{R(p)} - 1 = 
               R'(p) \int^p \frac 1{R(q)}\mbox dq  = R'(p) Q(p).
\end{eqnarray}
The time derivative of the internal energy follows from the chain rule:
\begin{eqnarray}
\paf{e_{int}}t  &=& 
e'_{int}(p) \paf p t = R'(p) Q(p) \paf pt =Q(p) \paf \rho t.
\end{eqnarray}
{\bf Time derivative of internal energy converted to spatial coordinates}\\[1ex]
Eliminating the time derivative using the continuity equation and using the symmetry property ${\sf DIVr}^*=-{\sf rGRAD}$, we find
\ContinuousDiscrete{
\begin{eqnarray}
\paf{e_{int}}t  &=&  -Q(p) \nabla \cdot \rho \vec v,
\end{eqnarray}\\[4ex]
}{
\begin{eqnarray}
\paf{{\tt e}_{int}}t  &=&  -\diag(Q({\tt p})) {\sf DIVr}~ {\tt v},
\\
\paf{E_{int}}t  &=&  \left<{\sf rGRAD}~ S({\tt p}), {\tt v}\right>_v.
\nonumber 
\end{eqnarray}
}
{\bf Energy equation}\\[1ex]
The time derivatives of local and total energies $e$ and $E$ are 
\ContinuousDiscrete{
\begin{eqnarray}
\paf et &=&
 -\nabla \cdot \frac \rho 2 |\vec v|^2 \vec v 
- \rho \vec v \cdot \nabla Q(p)
\nonumber \\ &&
-Q(p) \nabla \cdot \rho \vec v 
\nonumber \\ &=&
 -\nabla \cdot ( \frac 12 |\vec v|^2 +  Q(p))~\rho \vec v,
\end{eqnarray}
}{
\begin{eqnarray}
\paf Et &=&
-\left< {\tt v}, {\sf rGRAD} ~Q({\tt p}) \right>_v +
\nonumber \\ &=&
 \left<{\sf rGRAD}~ S({\tt p}), {\tt v}\right>_v
\nonumber \\ &=&0,
\end{eqnarray}
}
so the energy equation is 
\ContinuousDiscrete{
\begin{eqnarray}
\paf{e}t &+&
 \nabla \cdot ( \frac 12 |\vec v|^2 + \rho Q(p)) \rho \vec v  =0.
\end{eqnarray}
}{
\begin{eqnarray}
\paf Et =0.
\end{eqnarray}
}

\bibliographystyle{abbrv}
\bibliography{References}
\end{document}